\documentclass[11pt,fleqn]{article}
\usepackage{amsfonts, amsmath, epsfig, color, verbatim, url}
\usepackage[rightcaption]{sidecap}
\newcommand{\ol}{\setlength{\itemsep}{0pt.}\begin{enumerate}}
\newcommand{\eol}{\end{enumerate}\setlength{\itemsep}{-\parsep}}
\newcommand{\ignore}[1]{}
\title{On coset leader graphs of structured linear codes}
\author{Eran Iceland and Alex Samorodnitsky
\thanks{School of Engineering and Computer Science,
The Hebrew University of Jerusalem,
Jerusalem 91904, Israel.
Research partially supported by ISF
grants 1241/11 and 1724/15, and by BSF grant 2010451.}}

\begin{document}
\date{}
\maketitle


\newtheorem{THEOREM}{Theorem}[section]
\newenvironment{theorem}{\begin{THEOREM} \hspace{-.85em} {\bf :}
}%
                        {\end{THEOREM}}
\newtheorem{LEMMA}[THEOREM]{Lemma}
\newenvironment{lemma}{\begin{LEMMA} \hspace{-.85em} {\bf :} }%
                      {\end{LEMMA}}
\newenvironment{lemma1}{\begin{LEMMA}}
                      {\end{LEMMA}}
\newtheorem{COROLLARY}[THEOREM]{Corollary}
\newenvironment{corollary}{\begin{COROLLARY} \hspace{-.85em} {\bf
:} }%
                          {\end{COROLLARY}}
\newtheorem{PROPOSITION}[THEOREM]{Proposition}
\newenvironment{proposition}{\begin{PROPOSITION} \hspace{-.85em}
{\bf :} }%
                            {\end{PROPOSITION}}
\newtheorem{CLAIM}[THEOREM]{Claim}
\newenvironment{claim}{\begin{CLAIM} \hspace{-.85em}
{\bf :} }%
                            {\end{CLAIM}}
\newtheorem{DEFINITION}[THEOREM]{Definition}
\newenvironment{definition}{\begin{DEFINITION} \hspace{-.85em} {\bf
:} \rm}%
                            {\end{DEFINITION}}
\newenvironment{definition1}{\begin{DEFINITION} \rm}
                            {\end{DEFINITION}}
\newtheorem{EXAMPLE}[THEOREM]{Example}
\newenvironment{example}{\begin{EXAMPLE} \hspace{-.85em} {\bf :}
\rm}%
                            {\end{EXAMPLE}}
\newenvironment{example1}{\begin{EXAMPLE} \rm}
                            {\end{EXAMPLE}}
\newtheorem{CONJECTURE}[THEOREM]{Conjecture}
\newenvironment{conjecture}{\begin{CONJECTURE} \hspace{-.85em}
{\bf :} \rm}%
                            {\end{CONJECTURE}}
\newtheorem{MAINCONJECTURE}[THEOREM]{Main Conjecture}
\newenvironment{mainconjecture}{\begin{MAINCONJECTURE} \hspace{-.85em}
{\bf :} \rm}%
                            {\end{MAINCONJECTURE}}
\newtheorem{PROBLEM}[THEOREM]{Problem}
\newenvironment{problem}{\begin{PROBLEM} \hspace{-.85em} {\bf :}
\rm}%
                            {\end{PROBLEM}}
\newtheorem{QUESTION}[THEOREM]{Question}
\newenvironment{question}{\begin{QUESTION} \hspace{-.85em} {\bf :}
\rm}%
                            {\end{QUESTION}}
\newtheorem{REMARK}[THEOREM]{Remark}
\newenvironment{remark}{\begin{REMARK} \hspace{-.85em} {\bf :}
\rm}%
                            {\end{REMARK}}

\newtheorem{PROCEDURE}[THEOREM]{Procedure}
\newenvironment{procedure}{\begin{PROCEDURE} \hspace{-.85em} {\bf :}
\rm}%
                            {\end{PROCEDURE}}

\newcommand{\thm}{\begin{theorem}}
\newcommand{\lem}{\begin{lemma}}
\newcommand{\pro}{\begin{proposition}}
\newcommand{\clm}{\begin{claim}}
\newcommand{\dfn}{\begin{definition}}
\newcommand{\rem}{\begin{remark}}
\newcommand{\xam}{\begin{example}}
\newcommand{\cnj}{\begin{conjecture}}
\newcommand{\mcnj}{\begin{mainconjecture}}
\newcommand{\prb}{\begin{problem}}
\newcommand{\que}{\begin{question}}
\newcommand{\cor}{\begin{corollary}}
\newcommand{\prf}{\noindent{\bf Proof:} }
\newcommand{\ethm}{\end{theorem}}
\newcommand{\elem}{\end{lemma}}
\newcommand{\epro}{\end{proposition}}
\newcommand{\eclm}{\end{claim}}
\newcommand{\edfn}{\bbox\end{definition}}
\newcommand{\erem}{\bbox\end{remark}}
\newcommand{\exam}{\bbox\end{example}}
\newcommand{\ecnj}{\bbox\end{conjecture}}
\newcommand{\emcnj}{\bbox\end{mainconjecture}}
\newcommand{\eprb}{\bbox\end{problem}}
\newcommand{\eque}{\bbox\end{question}}
\newcommand{\ecor}{\end{corollary}}
\newcommand{\eprf}{\bbox}
\newcommand{\beqn}{\begin{equation}}
\newcommand{\eeqn}{\end{equation}}
\newcommand{\wbox}{\mbox{$\sqcap$\llap{$\sqcup$}}}
\newcommand{\bbox}{\vrule height7pt width4pt depth1pt}
\newcommand{\qed}{\bbox}
\newcommand{\overbar}[1]{\mkern 1.75mu\overline{\mkern-1.75mu#1\mkern-1.75mu}\mkern 1.75mu}
\def\sup{^}

\def\H{\{-1,1\}^n}
\def\B{\{0,1\}^n}

\def\S{S(n,w)}

\def\n{\lfloor \frac n2 \rfloor}

\def \E{\mathbb E}
\def \R{\mathbb R}
\def \Z{\mathbb Z}
\def \F{\mathbb F}
\def \T{\mathbb T}
\def \N{\mathbb N}

\def\<{\left<}
\def\>{\right>}
\def \({\left(}
\def \){\right)}
\def \e{\epsilon}
\def \r{\rfloor}

\def \F{{\mathbb{F}}}

\def \FF{{\cal F}}

\def \1{{\bf 1}}

\def \noi{\noindent}

\def\Tp{Tchebyshef polynomial}
\def\Tps{TchebysDeto be the maximafine $A(n,d)$ l size of a code with distance $d$hef polynomials}
\newcommand{\rarrow}{\rightarrow}

\newcommand{\larrow}{\leftarrow}

\overfullrule=0pt
\def\setof#1{\lbrace #1 \rbrace}

\begin{abstract}

We suggest a new approach to obtain bounds on locally correctable and some locally testable binary linear codes, by arguing that these codes (or their subcodes) have coset leader graphs with high discrete Ricci curvature.

The bounds we obtain for locally correctable codes are worse than the best known bounds obtained using quantum information theory, but are better than those obtained using other methods, such as the "usual" information theory. (We remark that our methods are completely elementary.)

The bounds we obtain for a family of locally testable codes improve the best known bounds.

\end{abstract}

\section{Introduction}

\noi We are interested in upper bounds on the cardinality of locally structured
linear subspaces of the Hamming space $\B$.

\noi To fix notions and some notation, let $C$ be a linear subspace of $\B$, and let $C^{\perp} =
\big\{y:~\<x,y\> = 0,~\forall x \in C\big\}$ be the dual space. We will assume that $C^{\perp}$
contains a rich family of vectors of constant length, and try to deduce that it is large
(alternatively, that $C$ is small).

\noi Specifically, we consider two families of locally
constrained linear binary codes. Such codes have numerous applications in theoretical computer
science (see \cite{dvir2014breaking} and the references therein). With that, essentially in all the cases, there is a significant gap between the
best known examples of such codes and upper bounds on their cardinality.


\noi Let $\FF$ be the family of local constraints on $C$ (that is, constant length
vectors in $C^{\perp}$). We will consider:

\begin{itemize}

\item {\it Locally correctable codes}

\noi These codes come with two parameters, an integer $q \ge 1$ and a density parameter $0 < \delta
< 1/q$. For each $1 \le i \le n$, the family $\FF$ contains either a unit vector supported at $i$, or at least $\delta n$ vectors of length
$q+1$ whose supports contain $i$ and are disjoint otherwise.

\item {\it Locally testable codes with high $3$-density}

\noi Here we have one integer parameter $\sigma$ which tends to infinity with $n$, and we assume that $\FF$
contains, for each $1 \le i \le n$, at least $\sigma$ vectors of length $3$ whose supports
contain $i$.

\end{itemize}

\subsection{Known bounds}

\subsubsection{Locally correctable codes}

\begin{itemize}
  \item For $q=1$, there are no locally correctable codes when $n$ is bigger than a constant  \cite{katz2000efficiency}.
  \item For $q=2$, the answer is $\Theta (\log n)$ \cite{goldreich2006lower}. \footnote{The constant hidden in the asymptotic notation here and below is allowed to depend on $\delta$.}
  \item For a constant $q>2$, there is a significant gap between upper and lower bounds.
  The best known locally correctable codes are the Reed-Muller codes of
  dimension $\Theta \( (\log n)^{q-1} \)$ \cite{macwilliams1977theory}.
  The best upper bound on the dimension is
  $O\left({n}^{\frac{\lceil{q/2}\rceil - 1}{\lceil{q/2}\rceil}}\right)$,
  up to a polylogarithmic factor
  \cite{kerenidis2004exponential, woodruff2007new, woodruff2012quadratic}.
\end{itemize}

\rem \label{rem:local-decod}

\noi In fact, all the upper bounds we mention hold also for {\it locally decodable codes}, which is yet another version of locally constrained codes.  For a locally decodable code $C$ of (unknown) dimension $D$, we assume that the first $D$ coordinates of a vector in $C$ determine the vector (in other words the
projection of a non-zero vector in $C$ on the first $D$ coordinates is non-zero). The parameters of
$C$ are $q$ and $\delta$, and the family of local constraints $\FF$ contains, for each $1 \le i \le
D$, at least $\delta n$ vectors of length $q+1$ whose supports contain $i$ and are disjoint
otherwise.

\noi Clearly, any locally correctable code is also locally decodable. The reverse implication does not hold \cite{kaufman2010locally}, but we
are not aware of any upper bounds separating these two families of binary codes.

\erem

\subsubsection{Locally testable codes}
Locally testable codes with high $3$-density were considered in \cite{ben2012towards}. Let us call a linear code $C$ {\it regular} if each column in its generating matrix appears with the same multiplicity. Then the following claim holds.

\begin{itemize}
\item
The dimension of a regular locally testable code with $3$-density $\sigma$ is at most
 $O\left(\frac{\log \sigma}{\sqrt \sigma} \cdot n\right)$ \cite{ben2012towards}.
\end{itemize}

\subsection{Our results}

\noi We follow the idea of Friedman and Tillich \cite{friedman2006generalized} and consider the {\it coset leader graph} of a binary linear code (see Definition~\ref{dfn:coset leader}). At this point it suffices to say that this is a Cayley graph whose cardinality is that of the corresponding code. Following a line of thought in \cite{friedman2006generalized} we view this graph as a homogeneous space and apply suitably modified tools from Riemannian geometry to upperbound its cardinality. Specifically, we will show this graph to have positive discrete Ricci curvature in the sense of \cite{ollivier2009ricci}. This will provide an upper bound on its diameter and hence on its size.

\noi {\bf Bounds for locally correctable codes:} In the statement of the next claim, and from now on, we will refer to a locally correctable code with parameters $q$ and $\delta$ as $q$-locally correctable (this in particular emphasizes the fact that $q$ is the more important parameter for the purpose of this discussion). Recall that we allow constants hidden in the asymptotic notation to depend on $\delta$.

\thm \label{thm:LCC}
Let $C$ be a $q$-locally correctable code with $q \geq 2$. Then the covering radius of
$C^\perp$ is $O\left(n ^{\frac{q-2}{q-1}}\right)$, and ${\rm dim} (C) \leq O\left(n
^{\frac{q-2}{q-1}}(\log n)^\frac{1}{q-1} \right)$.
\ethm

\noi While this is weaker than the best known bounds, we observe
that for $q > 3$ the exponent of $n$ in our bound on ${\rm dim}(C)$ lies strictly between that in
\cite{katz2000efficiency} obtained via Shannon's entropy, and the best known bound
\cite{kerenidis2004exponential,woodruff2007new}, which uses highly non-trivial facts, such as subadditivity of quantum entropy.

\noi For $q=3$ we do somewhat better. Theorem~\ref{thm:LCC} bounds the dimension of a $3$-locally correctable code by $\sqrt{n\log n}$, which is
only logarithmically weaker than the best known bound of $\sqrt n$ \cite{woodruff2012quadratic} (but stronger than $n^{2/3}$ of \cite{katz2000efficiency}). In fact, we can recover
the $\sqrt n$ bound in one special case.

\dfn \label{def:perfect_LCC} We say that a $q$-locally correctable code is \emph{perfect} if $\delta = \frac{n-1}{qn}$ for
every $i \in [n]$ (that is, the density parameter is as large as possible). \edfn

\thm \label{thm:3perLCC} Let $C$ be a perfect $3$-locally correctable code. Then ${\rm dim} (C) \leq O\(\sqrt
n\)$. \ethm

\noi {\bf Bounds for locally testable codes with high $3$-density:} We improve on the bounds for locally testable codes with high $3$-density.

\thm \label{thm:regLTC}
The dimension of a regular locally testable code with $3$-density $\sigma$ is at most
$\frac{2}{\sqrt{\sigma}} \cdot n$.
\ethm

\noi This bound is tight, up to a constant factor \cite{dinur2011dense}.

\noi We also consider a more general case in which the multiplicity of the columns is allowed to vary.

\thm \label{thm:bounded_repetitions} Let $v_1, v_2, \ldots , v_n$ be the columns of a generating
matrix $G$ of a code $C$. Let $p$ be the maximal multiplicity of a column in $G$. Assume that each
coordinate participates in at least $\sigma$ linear dependencies of length three, and that $\sigma > p$. Then
\[
\mbox{dim}(C) \le O\(\frac{\log\(\lceil\sigma/p\rceil\)}{\lceil\sigma/p\rceil} \cdot n\)
\]
\ethm

\noi This is tight, up to the $\log\(\lceil\sigma/p\rceil\)$-factor, see Example~\ref{ex:counter}.

\subsection{Our approach in more detail}

\noi Our starting point is the elegant proof of \cite{friedman2006generalized} for the first linear programming bound for binary linear codes. We start with the definition of coset leader graphs.

\dfn
\label{dfn:coset leader}
The {\it coset leader graph} $\T$ of a linear code $C \subseteq \B$ is the Cayley graph of the quotient group $\F_2^n / C^\perp$ with respect to the set of generators given by the standard basis $e_1+C^\perp, \ldots, e_n+C^\perp$.
\edfn
\noi Note that $\T$ has $|C|$ vertices. Note also that $\T$ may have loops or parallel edges (if $C^{\perp}$ contains non-zero vectors of Hamming weight less than $3$).

\noi \cite{friedman2006generalized} employs discrete versions of comparison theorems in Riemannian geometry comparing, on one hand, the growth of neighborhoods in $\T$ with the growth of neighborhoods in $\B$ and, on the other hand, the spectral behaviour of the Laplacian of $\T$ with the Laplacian of $\B$.

\noi Following \cite{friedman2006generalized}, we view $\T$ as a "discrete manifold", and try to estimate the cardinality of $\T$ by employing insights and tools borrowed from Riemannian geometry. The main technical notion we use is that of {\it discrete Ricci curvature}, due to Ollivier \cite{ollivier2009ricci}. We show that in the cases we consider $\T$ has 'high positive curvature' (much higher, say, than that of the Hamming cube)\footnote{To be more precise, in the proof of Theorem~\ref{thm:LCC} we consider the coset leader graph $\T$ of a {\it subcode} of $C$, rather than that of the whole code.}.

\noi This allows us to upper bound the diameter of $\T$ (equivalently, the covering radius of $C^\perp$), using a discrete version of the Bonnet-Myers theorem from Riemannian geometry \cite{ollivier2009ricci}. Since $\T$ is a regular graph of degree $n$, an upper bound on its diameter implies a bound on $|\T|$, and hence on $|C|$.

\rem
\label{rem:cmp}
While our approach is "curvature based", most of the bounds on local codes in the literature are based on isoperimetric inequalities or their information theoretic versions \cite{katz2000efficiency, kerenidis2004exponential, woodruff2007new}. In Riemannian geometry the notions of curvature and isoperimetry are closely related. Better isoperimetric inequalities for graphs with "positive discrete curvature" are known in the discrete setting as well \cite{lin2010ricci, bauer2015li, klartag2015discrete}. It seems natural to ask whether this connection might be exploited in order to improve coding bounds.
\erem

\section{The Main Technical Lemma}

\noi Our bounds are based on the following key lemma.

\begin{lemma1}
\label{lem:main}
Let $u_1,..., u_m$ and $v_1,..., v_n$ be correspondingly the rows and the columns of an $m \times n$ matrix over $\F_2$. Let
$V = {\rm span} (u_1,..., u_m)$.
Suppose that for each $i = 1,..., n$ with $v_i \not = 0$ there are at least $K>0$ disjoint pairs of indices $\{j,l\}$ such that $i \not \in \{j,l\}$ and $v_i = v_j + v_l$. Then, for the coset leader graph $\mathbb T = \{0,1\}^n / V^\perp$ of $V$ holds
\[
{\rm diam} (\T) \leq \frac{n}{K+1}.
\]
\end{lemma1}

\noi We observe that this implies a bound on the dimension of $V$.

\cor \label{cor:dimV}
$$\dim V \leq \log_2\(\sum_{i = 0}^{{\rm diam} (\T)} \binom{n}{i}\) \leq  \frac {n\log (K+1) + n/\ln 2}{K+1}$$.
\ecor
\prf (of Corollary~\ref{cor:dimV})

\noi Recall that $\T$ is an $n$-regular Cayley graph of an Abelian group. Hence
\[
|V| = |\T| \le \sum_{i = 0}^{{\rm diam} (\T)} \binom{n}{i}  \le 2^{n H\(\frac{1}{K+1}\)} \le 2^{\frac {n\log (K+1) + n/\ln 2}{K+1}}.
\]
Here $H(x) = x \log_2 \frac 1x + (1-x) \log_2 \frac{1}{1-x}$ is the binary entropy function. For the second inequality recall that for any $0 \le k \le r$ holds $\sum_{i = 0}^k \binom{r}{i} \le 2^{r H\(\frac kr\)}$ (Theorem 1.4.5. in \cite{van1999introduction}). For the third inequality, note that $(1-x) \ln \frac{1}{1-x} \le x$, for $0 \le x < 1$.
\eprf

\cor
\label{cor:q=2}
The case $q = 2$ of  Theorem~\ref{thm:LCC} holds.
\ecor
\prf (of Corollary~\ref{cor:q=2})

\noi Let $\T$ be the coset leader graph of a locally correctable code $C$ with parameters $2$ and $\delta$. By definition, a generating matrix of $C$ satisfies the assumptions of Lemma~\ref{lem:main} with $K = \delta n$, and hence ${\rm diam}(\T) \le 1/\delta$ and ${\rm dim}(C) \le O\(\frac{\log n}{\delta}\)$. Since it is easy to see that the diameter of $\T$ is precisely the covering radius of $C^\perp$, this proves the case $q=2$ of Theorem~\ref{thm:LCC}.
\eprf

\subsection{Examples}

\xam
\label{ex:hadamard}
Let $m$ be a positive integer and let $n = 2^m - 1$. The generating matrix of the {\it Hadamard code} $C$ of length $n$ \cite{macwilliams1977theory} is the $m \times n$ matrix whose columns are all the non-zero vectors in $\F_2^m$. Let $\T$ be the coset leader graph of $C$. Since the columns of the generating matrix are non-zero and distinct, $\T$ is a simple $n$-regular graph with $2^m = n+1$ vertices, namely it is the complete graph on $n+1$ vertices.
In this case the assumptions of Lemma~\ref{lem:main} hold with $K = (n-1)/2$, and it gives the tight bounds ${\rm diam}(\T) \le \lfloor \frac{2n}{n+1} \rfloor = 1$ and $|\T| \le n+1$.   \exam

\xam
\label{ex:hadamard-times-2}
Let $C$ be the direct product of two Hadamard codes. That is, assume $n = 2 \cdot \(2^m  - 1\)$, and let the generating matrix of $C$ be a $2m \times n$ block-diagonal matrix with two $m \times (n/2)$ blocks whose columns are all the non-zero vectors in $\F_2^m$. In this case $\T$ is the Cartesian product of two complete graphs on $n/2+1$ vertices. That is, $|\T| = \frac{(n+2)^2}{4}$ and ${\rm diam}(\T) = 2$. The conditions of the lemma hold with $K = (n-2)/4$, leading to an upper bound of $3$ on the diameter and of $\binom{n}{3} + \binom{n}{2} + n + 1$ on the cardinality of $\T$.
\exam

\noi In the remainder of this section we proceed as follows. We start with comparing Lemma~\ref{lem:main} to related results in the literature. Next, we describe the key notion of discrete curvature on graphs. Finally, we prove the lemma in Section~\ref{sec:prf_main}.

\subsection{Lemma~\ref{lem:main} and related results}

\noi In this subsection we expand on Remark~\ref{rem:cmp}. Upper bounds on locally correctable codes in the literature follow from bounds on locally decodable codes. (Our
approach applies directly to locally correctable codes, which might explain its relative simplicity.) A typical approach uses isoperimetric inequalities. \cite{goldreich2006lower, woodruff2007new, woodruff2012quadratic} use a weighted version of the edge-isoperimetric inequality on the boolean cube (\cite{goldreich2006lower}). Another version of the edge-isoperimetric inequality is proved and used in \cite{ben2012towards}. We compare Lemma~\ref{lem:main} and
Corollary~\ref{cor:dimV} with these two results, which we restate in our language.

\begin{lemma1} {\normalfont (\cite[Lemma 3.3]{goldreich2006lower}){\bf:}}
\label{lem:combinatorial}
Let $u_1,..., u_m$ and $v_1,..., v_n$ be, correspondingly, the rows and the columns of an $m \times n$ matrix over $\F_2$. Let
$V = {\rm span} (u_1,..., u_m)$ and assume that $v_1,..., v_{{\rm dim} V}$ span the column space.
Suppose that for each $i = 1,..., \dim V$ there are $K_i>0$ disjoint pairs of indices $\{j,l\}$ such that $i \not \in \{j,l\}$ and $v_i = v_j + v_l$. Let $K = \(\sum_{i=1}^{\dim V}K_i\) / \dim V$. Then
\[
\dim V \leq \frac {n \log n}{2K}.
\]
\end{lemma1}

\begin{lemma1} {\normalfont (\cite[Lemma 3.15]{ben2012towards}){\bf:}}
\label{lem:BV}
Let $u_1,..., u_m$ and $v_1,..., v_n$ be, correspondingly, the rows and the columns of an $m \times n$ matrix over $\F_2$. Let
$V = {\rm span} (u_1,..., u_m)$ and assume that $v_1,..., v_{{\rm dim} V}$ span the column space.
Suppose that for each $i = 1,..., \dim V$ there are $K > 0$ disjoint pairs of indices $\{j,l\}$ such that $i \not \in \{j,l\}$ and $v_i = v_j + v_l$. Then
\[
\dim V \leq \frac {n \log K +n}{K}.
\]
\end{lemma1}

\noi We collect the assumptions and the conclusions of the three claims in the following table (omitting constants for readability):

\vspace{1em}
\begin{center}
\begin{tabular}{|l||c|c|c|}
  \hline
   &Lemma~\ref{lem:main} + Cor.~\ref{cor:dimV}& Lemma~\ref{lem:combinatorial} & Lemma~\ref{lem:BV} \\ \hline
  Disjoint rep's for: & all columns & basis & basis \\
  At least $K$ rep's for: & all columns & average basis column & all basis columns  \\
  Dimension at most & $\frac {n \log K}{K} $& $\frac{n \log n}{K}$ & $\frac {n \log K}{K}$ \\
  Diameter at most & $\frac nK$ & $\frac{n \log n}{K}$ & $\frac {n \log K}{K}$ \\
  \hline
\end{tabular}
\end{center}

\vspace{1em}

\noi This table requires some reading help, which we provide here. The first two rows present the assumptions, and the last two the bounds. In this context, a family of disjoint representations of a column $i$ is a collection of disjoint pairs of indices $\{j,l\}$ such that $i \not \in \{j,l\}$ and $v_i = v_j + v_l$.

\noi The first row specifies whether such family is assumed to exist for all column vectors $v_1, \ldots, v_n$ or only for the basis $v_1, \ldots, v_{\dim V}$. The second row indicates whether the lower bound $K$ is on the minimal or the average size of a family (over the relevant coordinates).

\noi The third row bounds the dimension of $V$. The fourth row bounds the diameter of the coset leader graph $\T = \B / V^\perp$, which is the same as the covering radius of $V$. Since lemmas~\ref{lem:combinatorial} and~\ref{lem:BV} do not consider the covering radius, we have filled out the corresponding entries using the fact that the covering radius of a linear code is upper bounded by its dimension \cite[Theorem 2.1.9]{cohen1997covering}.

\noi The next example shows that in Lemma~\ref{lem:main} the assumption on minimal family size cannot be replaced by that on average family size, without affecting the bounds.

\xam
For an integer $m$, let $n = 2^m -1$, and let $A$ be the generating matrix of the Hadamard code (Example~\ref{ex:hadamard}) with $m$ rows and $n$ columns. Let $I_m$ be the $m \times m$ identity matrix. Consider a linear code defined by the following $(2m) \times (n+m)$ generating matrix:
$$
\(
\begin{array} {ccc|c}
   & A &  & 0  \\ \hline
   & 0 &  & I_m \\
\end{array}
\)
$$

\noi In this case, the average family size is linear in $n$ (but the minimal family size is zero). The coset leader graph is the Cartesian product of the complete graph on $n+1$ vertices with the $m$-dimensional discrete cube. Hence its diameter is logarithmic in $n$ (as opposed to constant).
\exam

\subsection{Discrete Curvature on Graphs}
\label{sec:Olliv}

\noi

\noi There are several possible ways to extend the notion of Ricci curvature from Riemannian geometry to the general setting of metric spaces and, in particular, graphs \cite{charney1994metric, ollivier2009ricci, petrunin2011alexandrov, bauer2012ollivier}. We use the approach of Ollivier \cite{ollivier2009ricci, ollivier2012curved}.
In the following discussion $G$ is a finite multigraph with a probability measure $m_x$ on its vertex set $V = V(G)$ assigned to each vertex $x$. We denote by $d$ the graph (shortest path) metric defined by $G$ on $V$.

\noi Recall that the {\it transportation distance} between two probability measures $\mu$ and $\nu$ on $V$ is defined as
\begin{equation}
\label{eq:td}
W_1 \(\mu, \nu\) = \min_q \sum_{\(x',y'\) \in G \times G} q(x',y')d(x',y')
\end{equation}
where the minimum is taken over all probability measures $q$ on the product space $V \times V$ whose marginals are $\mu$ and $\nu$.

\dfn
Let $x \not = y \in V$. The \emph{coarse Ricci curvature} $\kappa(x,y)$ along $(x,y)$ is
\beqn
\label{eq:kappa}
\kappa(x,y) = 1 - \frac{W_1\(m_x,m_y\)}{d(x,y)}
\eeqn
\edfn

\noi A canonical choice for the measure $m_x$ is the uniform probability measure on the metric ball of radius $1$ around $x$. In this case the coarse Ricci curvature $\kappa(x,y)$ along $(x,y)$ is positive if and only if the mean distance between the metric balls around $x$ and $y$ (as measured by $W_1 (m_x, m_y)$) is smaller than the distance between $x$ and $y$. This conforms to the intuition that in spaces with positive curvature metric balls are closer on average than their centers (and vice versa for spaces with negative curvature).

\noi The {\it curvature} $\kappa(G)$ of the {\it graph} $G$ is defined as the minimum of $\kappa(x,y)$ over all pairs of vertices. This minimum is attained on a pair of adjacent vertices \cite[Proposition 19]{ollivier2009ricci}.\footnote{This, and Proposition~\ref{prop:BM} below are simple consequences of the triangle inequality for the transportation distance.} In particular, curvature is a local property.

\noi The key claim we need is the following discrete version of the Bonnet-Myers theorem for Riemannian manifolds\footnote{The classical Bonnet-Myers theorem for Riemannian manifolds states that if the Ricci curvature of an $n$-dimensional complete Riemannian manifold $M$ is at least $(n - 1)\kappa > 0$, then the manifold is compact and its diameter is at most $\pi/\sqrt \kappa$.}  \cite[Proposition 23]{ollivier2009ricci}. For $x \in V$, let $\delta_x$ be the probability measure concentrated on $x$, and let $J(x) = W_1(\delta_x, m_x)$. (E.g., if $G$ is a simple $n$-regular graph, and $m_x$ is the uniform measure on the metric ball of radius $1$ around $x$, then $J(x) = \frac{n}{n+1}$.)

\pro {\normalfont (\cite[Proposition 23]{ollivier2009ricci})}
\label{prop:BM}
For any graph $G$ and a family of probability measures $\{\mu_x\}_{x \in V(G)}$ holds
\[
{\rm diam}(G) \le \frac{2 \cdot \max_{x \in G} J(x)}{\kappa(G)}.
\]
\epro

\rem
Many of these ideas appear also in the theory of random walks on graphs, see \cite{levin2009markov}, especially chapter 14.
\erem

\subsection{Proof of Lemma~\ref{lem:main}}
\label{sec:prf_main}

\noi In this section we prove Lemma~\ref{lem:main}. We choose a family of probability measures $\{\mu_x\}_{x \in U}$ on the vertex set $U$ of $\T$ which enables us to bound the coarse Ricci curvature on $\T$ from below, and then apply Proposition~\ref{prop:BM}.

\noi Recall that $\T$ is an $n$-regular multigraph. For $x \in U$, we define $m_x$ to be the measure induced by the uniform measure on the edges incident to $x$. That is, for $y$ adjacent to $x$ we set $m_x(y)$ to be the number of edges between $x$ and $y$, divided by $n+1$; and we let $m_x(x)$ be the number of loops at $x$ {\it plus one}, divided by $n+1$. The measure $m_x$ is supported on the metric ball of radius $1$ around $x$ and, if $\T$ is a simple graph, then $m_x$ is uniform on this set.

\noi We observe that the local structure of $\T$ at (any) vertex $x$, and hence the measure $m_x$, can be described in terms of the column vectors $v_1, \ldots, v_n$. In fact, the number of loops at any vertex of $\T$ equals to the number of zero vectors among $v_1, \ldots, v_n$. Similarly, the number of edges between two distinct vertices $x$ and $x + e_i$ (the addition is in the factor group $\F^n_2 / V^{\perp}$) is the number of times $v_i$ appears as a column vector.

\noi Next, we upper bound the transportation distance between measures $m_x$ and $m_y$, for distinct adjacent vertices $x$ and  $y$. Let $y = x+e_i$ for some $1 \leq i \leq n$. By the assumption of the lemma, there are some $M \ge K$ disjoint pairs of indices $\{j,l\} \subseteq [n] \setminus \{i\}$ such that $v_j + v_l = v_i$. Equivalently, $e_i + e_j + e_l \in V^{\perp}$ and hence
$$x + e_j = y + e_l \quad \quad {\rm and} \quad \quad x+e_l = y+ e_j$$

\noi This means that the points $x+e_j$ and $x+e_l$ belong to the supports of both $m_x$ and $m_y$, and we have identified an overlap between the two measures, of weight $\frac{2}{n+1}$ in each measure. Going over all the $M$ representations $v_j + v_l = v_i$ produces an overlap of weight $\frac{2M}{n+1}$ in each of the measures.

\noi The identity $y = x+e_i$ gives an additional overlap of at least $\frac{2}{n+1}$ between the measures. This brings the total overlap to at least $\frac{2M+2}{n+1}$.

\noi We now transport $m_x$ to $m_y$ as follows. The points in the joint support stay in place. All the remaining mass in $m_x$ is moved by a unit distance in parallel. That is, we move a point $z$ in the unit ball around $x$ to the point $z + e_i$ in the unit ball around $y$.

\noi Computing the total amount of work gives
\[
W_1(m_x, m_y) \leq \frac{(n+1) - (2M+2)}{n+1} \leq \frac{(n+1) - (2K+2)}{n+1}.
\]

\noi Hence, by (\ref{eq:kappa}), the coarse Ricci curvature along $(x,y)$ is at least $\frac{2(K+1)}{n+1}$. Since this holds for any adjacent pair of vertices, we have $\kappa(\T) \ge \frac {2(K+1)}{n+1}$.

\noi Applying Proposition~\ref{prop:BM} (note that $J(x) \leq \frac{n}{n+1}$ for all $x \in U$) gives
\[
\mbox {diam} \(\T\) \leq \frac {\frac{2n}{n+1}}{\kappa(\T)} \le \frac{n}{K+1},
\]
concluding the proof of the lemma.
\eprf

\section{Bounds on Locally Correctable Codes}
\label{sec:prf_LCC}

\subsection{Proof of Theorem~\ref{thm:LCC}}

\noi The case $q = 2$ of the theorem is treated in Corollary~\ref{cor:q=2}. In this section we deal with larger values of $q$.

\noi Let $C$ be a locally correctable code with parameters $q > 2$ and $\delta$. Fix a generating matrix of $C$ and let its columns be $v_1, \ldots ,v_n$. Let $N = \{i \mid v_i \not = 0\}$. By definition, for each coordinate $i \in N$, there is a family $M_i$ of at least $\delta n$ disjoint $q$-subsets of $[n] \setminus \{i\}$ such that the vectors indexed by each subset sum to $v_i$.

\noi Our argument works (essentially) by reduction to the base case $q = 2$. Let us start with a quick overview. We will show that there is a subset $B$ of $[n]$ such that for any $i \in N \setminus B$ there are many $q$-tuples $\alpha \in M_i$ with $|\alpha \setminus B| \leq 2$. Contracting by $B$ will produce a code whose generating matrix satisfies the conditions of Lemma~\ref{lem:main}, with a parameter $K$ related to the parameters of the original code. Applying Lemma~\ref{lem:main} and Corollary~\ref{cor:dimV} will complete the proof. Let us mention that this approach is similar to that in \cite{woodruff2007new} and \cite{dinur2011dense}.

\lem \label{lem:cl:B}
Let $q>2$. For each $1 \leq a \leq (\log n)^{1/(q-1)}$ there exists a subset $B \subseteq [n]$, such that:
\begin{itemize}
  \item $|B| \quad \leq \quad  \Big( a + \frac{4}{\delta a^{q-2}} \Big) \cdot n^{\frac{q-2}{q-1}}$
  \item For every $i \in N \setminus B$ holds  $~~~\bigg|\Big\{\alpha \in M_i \ : \ |\alpha \setminus B| \leq 2 \Big\}\bigg| \quad \geq \quad \frac{\delta}{2} a^{q-2} \cdot n^{\frac{1}{q-1}}$
\end{itemize}
\elem

\prf

\noi Set $\theta = a \cdot n^{-\frac{1}{q-1}}$ and observe that $0 < \theta < 1$. We construct a random subset $B \subseteq [n]$ satisfying the assertions of the lemma in two steps. In the first step we add to $B$ elements in $[n]$ chosen independently at random with probability $\theta$. With high probability, this will produce a set of cardinality about $n \theta = a \cdot   n^{\frac{q-2}{q-1}}$, satisfying the second claim of the lemma for all but a small number of indices $i \in N$. In the second step we will add to $B$ all these exceptional indices and in this way ensure that both claims of the lemma hold.

\noi Let $X_1, \ldots, X_n$ be i.i.d. Bernoulli random variables,
$X_j = \left\{
  \begin{array}{ll}
    1, & \hbox{w.p. $\qquad \theta$;} \\
    0, & \hbox{w.p. $\quad 1-\theta$.}
  \end{array}
\right.
$

\noi Let $B_0 = \{ j \: : \: X_j = 1\}$. For a $q$-subset $\alpha \subseteq [n]$, let $W_{\alpha}$ be indicator of the event $|B_0 \cap \alpha| \ge q-2$.

\noi For $i \in N$, let $Y_i = \sum_{\alpha \in M_i} W_{\alpha}$. Since the $q$-tuples in $M_i$ are disjoint, the random variables $\{W_{\alpha}\}_{\alpha \in M_i}$ are independent, and hence $Y_i$ is a binomial random variable with parameters $|M_i|$ and $\eta = Pr\(W_{\alpha} = 1\) > \theta^{q-2}$. In particular,
$\E(Y_i) = |M_i|\cdot \eta > \delta n \theta^{q-2} = \delta a^{q-2} \cdot n^{\frac{1}{q-1}}$.

\noi Hence, by Chebyshev's inequality,
\[
  \Pr \(Y_i < \frac{\delta}{2} a^{q-2} \cdot n^{\frac{1}{q-1}}\) \leq
  \Pr \(Y_i < \frac{{\E}(Y_i)}{2} \)      \leq
  \frac {{\rm var} (Y_i)}{({\E}(Y_i)/2)^2}\leq
  \frac {4}{{\E} (Y_i)}         \leq
  \frac{4}{\delta a^{q-2} \cdot n^{\frac{1}{q-1}}}
\]

\noi Let $B = B_0 \bigcup \left\{i \in N \: : \: Y_i < \frac{\delta}{2} a^{q-2} \cdot n^{\frac{1}{q-1}} \right\}$.

\noi By the definition of $B$, for all $i \in N \setminus B$ holds
\[
\bigg|\Big\{\alpha \in M_i \:: \: |\alpha \setminus B| \leq 2 \Big\}\bigg| \geq
\bigg|\Big\{\alpha \in M_i \:: \: |\alpha \setminus B_0| \leq 2 \Big\}\bigg|  = Y_i \geq
\frac{\delta}{2} a^{q-2} \cdot n^{\frac{1}{q-1}}
\]

\noi Therefore, $B$ satisfies the second claim of the lemma. To verify that for some choice of $B$ the first claim holds as well, we upperbound the expectation of $|B|$ appropriately.
\[
  \E(|B|) = \E(|B_0|) + \E (|B\setminus B_0|) \le
    n \theta + \sum_{i \in N} \Pr \( Y_i < \frac{\delta}{2} a^{q-2} \cdot n^{\frac{1}{q-1}}\) \leq \( a + \frac{4}{\delta a^{q-2} }\) \cdot n^{\frac{q-2}{q-1}}
\]
\eprf

\noi We proceed with the proof of Theorem~\ref{thm:LCC}. Let $a$ be a parameter in the interval$\left[1,(\log n)^{1/(q-1)}\right]$ (we will optimize over the value of $a$ later on). Let $B = B(a) \subseteq [n]$ be the subset of indices given by Lemma~\ref{lem:cl:B}. Let $U = U(B) = {\rm Span} (\{v_i \,:\, i \in B\})$. Let $C_B$ be the subcode of $C$ containing  the vectors in $C$ which vanish on $B$. Let $\T = \F^n_2 / C^{\perp}$ and $\T_B = \F^n_2 / C_B^{\perp}$ be the coset leader graphs of $C$ and $C_B$ respectively. Then the following holds.

\lem
\label{lem:factor-code}
\begin{itemize}

\item
Any generating matrix of $C_B$ satisfies the conditions of Lemma~\ref{lem:main} with $K = \frac{\delta}{2} a^{q-2} \cdot n^{\frac{1}{q-1}}$.

\item
$\dim C = \dim C_B + \dim U$.

\item
${\rm diam} (\T) \leq {\rm diam} (\T_B) + \dim U$.

\end{itemize}
\elem
\prf
We start with the first claim. Let $G_B$ be a generating matrix of $C_B$. Note that the preceding discussion, and in particular the choice of the set $B$, has been independent of the generating matrix of $C$ we have chosen, and hence we may assume that $G_B$ is a row submatrix of this generating matrix, which we will denote by $G$. Let $u_1,...,u_n$ be the columns of $G_B$. Then $u_i$ is a restriction of $v_i$ to a (fixed) subset of coordinates for all $1 \le i \le n$.

\noi Let $u_i$ be a non-zero column of $G_B$. We need to show that there are at least $K$ disjoint pairs of indices $\{j,l\}$ with $u_i = u_j + u_l$. First, note that $i \in N \setminus B$. Indeed, $u_i$ is zero for $i \in B$, by the definition of $C_B$, and $v_i$ (and hence $u_i$) is zero for $i \not \in N$, by the definition of $N$.

\noi Since $i \in N \setminus B$, there are at least $K$ disjoint $q$-tuples $\alpha \in M_i$ with $|\alpha \setminus B| \leq 2$. We will find a coordinate pair $\{j,l\}$ with $u_i = u_j + u_l$ contained in each of these tuples, and this will complete the argument. Fix $\alpha$. By definition, $\sum_{k \in \alpha} v_k = v_i$, implying $\sum_{k \in \alpha} u_k = u_i$. Since $u_s = 0$ for $s \in B$, this means $\sum_{k \in \alpha \setminus B} u_k = u_i$. Since $u_i \not = 0$, the set $\alpha \setminus B$ is not empty. If $|\alpha \setminus B| = 2$, take $\{j,l\} = \alpha \setminus B$. If $|\alpha \setminus B| = 1$, take $j$ to be the unique element of $\alpha \setminus B$, and $l$ any element of $\alpha \cap B$.

\noi The second claim is a well-known fact in linear algebra. We provide a brief argument for completeness. Right multiplication by $G$ defines an isomorphism between $\F_2^{\dim C}$ and $C$. The claim is implied by the observation that the pre-image of $C_B$ under this isomorphism is precisely $U^{\perp}$.

\noi We pass to the third claim. Since the diameter of the coset graph of a code equals to the covering radius of the dual code, the claim is that the covering radius of $C^{\perp}$ is upper bounded by the covering radius of $C_B^{\perp}$ plus the dimension of $U$. We will show this by finding, for each vector $x \in C_B^{\perp}$, a vector $y \in C^{\perp}$ such that $|x - y| \le \dim U$. Observe that $C_B^{\perp} = \Big\{x \in \B, ~\sum_{i=1}^n x_i v_i \in U\Big\}$. Let $x \in C_B^{\perp}$, and let $\sum_{i=1}^n x_i v_i = u \in U$. The vector $u$ can be written as a linear combination of columns in $B$, of length at most $\dim U$. Let $z \in \B$ be the characteristic vector of this linear combination. Then $|z| \le \dim U$ and $y = x + z \in C^{\perp}$, completing the proof.

\eprf

\noi Now we are ready to complete the proof of Theorem~\ref{thm:LCC}. To bound the covering radius of $C^\perp$, which is the same as the diameter of $\T$, take $a = 1$. This gives $|B| = \(1 + \frac{4}{\delta}\) \cdot n^{\frac{q-2}{q-1}}$ in Lemma~\ref{lem:cl:B} and $K = \frac{\delta}{2} \cdot n^{\frac{1}{q-1}}$ in Lemma~\ref{lem:factor-code}. By Lemma~\ref{lem:main}, ${\rm diam} (\T_B) \leq \frac{n}{K+1}$, and hence
\[
{\rm diam} (\T) \le {\rm diam} \T_B + \dim U  \le \frac{n}{K+1} + |B| \le  O\(n^{\frac{q-2}{q-1}}\).
\]

\noi To bound the dimension of $C$, take $a = (\log n)^\frac{1}{q-1}$. This gives $|B| \approx  n^{\frac{q-2}{q-1}} (\log n)^\frac{1}{q-1}$ and $K = \frac {\delta}{2} n^{\frac {1}{q-1}} (\log n)^\frac{q-2}{q-1}$. By Corollary~\ref{cor:dimV}, $\dim C_B \leq O \(  n^{\frac{q-2}{q-1}} (\log n)^{\frac {1}{q-1}}\)$, and hence,
\[
\dim C \le \dim C_B + \dim U \le \dim C_B + |B| \le  O \(  n^{\frac{q-2}{q-1}} (\log n)^{\frac {1}{q-1}}\).
\]

\subsection{Proof of Theorem~\ref{thm:3perLCC}}\label{subsec:prf_3perLCC}

\noi Let $\T = \B / C^\perp$ be the coset leader graph of $C$. We will show that the neighborhoods of (any) vertex in $\T$ grow rather slowly, which will imply that $\T$, and hence $C$, are not too large. For $r \ge 0$, let $S_r^{\T}$ be the sphere of radius $r$ around $C^\perp$ in $\T$. The key observation is that there are many edges in $\T$ between the consecutive spheres $S_{r-1}^{\T}$ and $S_{r}^{\T}$.


\lem \label{lem:bipartite}
Let $r \ge 2$. Assume that $S_r^{\T}$ is not empty. Then there are at least $\(\lfloor r/2 \rfloor\)^2$ edges between any vertex $x+C^\perp \in S_r^{\T}$ and $S_{r-1}^{\T}$.
\elem

\rem
\label{rem:compare}
This should be compared to the situation in the discrete cube $\B$, also an $n$-regular graph, in which a vertex at distance $r$ from zero is connected to the sphere of radius $r-1$ around zero by exactly $r$ edges.
\erem

\noi Before proving the lemma, let us show that it implies the claim of the theorem. By the lemma, there are at least $\(\lfloor r/2 \rfloor\)^2 \cdot |S_r^{\T}|$ edges between $S_{r-1}^{\T}$ and $S_{r}^{\T}$. On the other hand, $\T$ is an $n$-regular graph, which means that there are at most $n \cdot |S_{r-1}^{\T}|$ such edges. Hence $\(\lfloor r/2 \rfloor\)^2 \cdot |S_r^{\T}| \leq n \cdot |S_{r-1}^{\T}|$, and this holds for any $r \ge 2$.

\noi The sphere of radius $1$ is of cardinality at most $n$. Multiplying consecutive inequalities provides an upper bound on the cardinality of a sphere of radius $r \ge 2$:
\[
\Big | S_r^{\T} \Big | ~~\leq~~ n \cdot \prod_{t = 2}^r \frac{n}{\(\lfloor t/2 \rfloor\)^2} ~~\leq~~ \(\frac{cn}{r^2}\)^r
\]
for an appropriate constant $c > 0$. The second inequality can be deduced e.g., from Stirling's formula. It is easy to see that this implies $|\T| = \sum_r \Big | S_r^{\T} \Big | \le c^{\sqrt{n}}$, for a (possibly different) constant $c$, completing the proof of the theorem.

\prf (of Lemma~\ref{lem:bipartite}).

\noi By assumption, $C$ is a perfect $3$-locally correctable code. This means that $n$ is $1$ modulo $3$, and that for all $1 \le i \le n$ there is a family $M_i$ of $\frac{n-1}{3}$ disjoint $3$-tuples partitioning $[n] \setminus \{i\}$, so that for any such $3$-tuple $\alpha$ holds $e_i + \sum_{j \in \alpha} e_j \in C^\perp$. That is, for any two indices $i < j$ there is a unique pair of indices $k \not = l$ such that $(j,k,l) \in M_i$. In particular, $\{i,j\} \cap \{k,l\} = \emptyset$, and $e_i + e_j + e_k + e_l \in C^\perp$.

\noi Let now $x + C^\perp \in S_{r}^{\T}$. We may assume that $x$ is of minimal weight in its coset, meaning that the Hamming weight of $x$ is $r$. We will also assume, for simplicity, that $x_i = 1$ for $1 \le i \le r$ (and $x_i = 0$ for $i>r$).

\noi For two indices $i,j$ with $1 \leq i < j \leq r$, let $k,l$ be such that $(j,k,l) \in M_i$. Note that this necessarily means that $x_k = x_l = 0$ (that is $k, l > r$). Indeed, otherwise $x' = x + e_i + e_j + e_k + e_l$ would be a vector in $x + C^\perp$ of weight smaller than $r$.
The key point for us is that the edges from $x+C^\perp$ in the directions $k,l $ lead down to $S_{r-1}^{\T}$. In fact, the vector $x + e_k = x + e_i + e_j + e_l$ is of weight $r-1$ (and similarly for $x + e_l$).

\noi Let $V \subseteq [n]$ contain all directions leading from $x$ down to $S_{r-1}^{\T}$. As we have seen, each pair of indices $i,j$ with $1 \leq i < j \leq r$ defines a pair $(k, l) \in V \times V$, which we interpret as an edge with vertices in $V$. From now on we assume, for simplicity, that $r$ is even. Going over $i,j$ with $1 \le i \le r/2 < j \le r$ defines a multigraph $G$ on $V$ with $r^2/4$ edges. In fact, we claim that $G$ is a simple graph, that is distinct pairs $i,j$ and $i_1,j_1$ define distinct edges $(k,l)$ and $(k_1,l_1)$. Indeed, otherwise $e_i + e_j + e_{i_1} + e_{j_1} \in C^\perp$, which means that $x' = x + e_i + e_j + e_{i_1} + e_{j_1}$ is a vector in $x + C^\perp$ of weight smaller than $r$.

\noi Next, we claim that $G$ is a disjoint union of stars. This would mean that the number of vertices of $G$ is larger than its number of edges, i.e., $|V| > r^2/4$, proving Lemma~\ref{lem:bipartite}. This claim is a simple corollary of the following auxiliary lemma.

\lem
\label{lem:stars}
Any edge of $G$ contains a vertex of degree $1$.
\elem

\prf
Assume to the contrary that there exists an edge $(k,l)$ in $G$ such that both $k$ and $l$ have degree at least $2$. There are two possible cases. Either $G$ contains a simple path $k_1 \rarrow k \rarrow l \rarrow l_1$ of length $4$, or $G$ contains a triangle with vertices $k, l, m$. Consider the first case. Let $(i,j)$, $1 \leq i \le r/2 < j \leq r$, be the pair of indices defining the first edge of the path, let $(i_1,j_1)$ define the second edge, and $(i_2,j_2)$ the third edge. We claim that $i \not = i_1$. Indeed, otherwise both $(j,k_1,k)$ and $(j_1,k,l)$ would be in $M_i$, contradicting the fact that $M_i$ is a family of disjoint triples. Next, we claim that $j = j_1$. If not, we would have
\[
\Big(e_i + e_{i_1} + e_j + e_{j_1}\Big) + \Big(e_{k_1} + e_l\Big) ~~ = ~~ \Big(e_i + e_j + e_{k_1} + e_k\Big) + \Big(e_{i_1} + e_{j_1} + e_{k} + e_l\Big) \in C^\perp,
\]
which would give us a vector $x' = x + \Big(e_i + e_{i_1} + e_j + e_{j_1}\Big) + \Big(e_{k_1} + e_l\Big)$ in $x + C^\perp$ of weight smaller than $r$.

\noi A similar argument shows that $i_1 \not = i_2$ and $j_1 = j_2$ (and therefore also $j = j_2$). We now observe that $i$ and $i_2$ also have to be distinct. Indeed, otherwise we would have both $(j,k_1,k)$ and $(j,l,l_1)$ in $M_i$.

\noi Taking everything into account, this means that
\[
\Big(e_i + e_{i_1} + e_{i_2} + e_j\Big) + \Big(e_{k_1} + e_{l_1}\Big) ~~ = ~~ \Big(e_i + e_j + e_{k_1} + e_k\Big) +...+ \Big(e_{i_2} + e_{j_2} + e_{l} + e_{l_1}\Big) \in C^\perp,
\]
giving a vector $x' = x + \Big(e_i + e_{i_1} + e_{i_2} + e_j\Big) + \Big(e_{k_1} + e_{l_1}\Big)$ in $x + C^\perp$ of weight smaller than $r$, and in this way reaching a contradiction.

\noi The second case of the lemma is similar (but simpler). We omit the analysis. This completes the proof of Lemma~\ref{lem:stars} and of Lemma~\ref{lem:bipartite}.

\eprf

\section{Bounds on Locally Testable Codes}

\noi In this section we prove Theorems~\ref{thm:regLTC}~and~\ref{thm:bounded_repetitions}. The proofs of both theorems are based on the following lemma.

\lem
\label{lem:bounded_repetitions}
Let $G$ be a matrix satisfying the assumptions of Theorem~\ref{thm:bounded_repetitions}. Then $G$ satisfies the assumptions of Lemma~\ref{lem:main} with $K=\lceil\sigma/p\rceil$.
\elem

\prf
\noi Let $t$ be the number of distinct columns of $G$ and assume, without loss of generality, that $v_1, \ldots ,v_t$ are pairwise distinct. That is, the first $t$ columns represent all the distinct columns in $G$. For $1 \leq i \leq t$, let $w_i$ denote the multiplicity of $v_i$ in $G$. Note that $\sum _{i=1}^t w_i = n$. We may, and will, assume that $w_1 \leq \ldots \leq w_t = p$.
For $1 \leq i \leq t$ with $v_i \not = 0$, let $N_i = \{(j,k):~1 \le j < k \le t,~  v_i = v_j + v_k\}$.

\noi Fix an index $1 \leq i \leq n$ with $v_i \not = 0$. We need to show that there are at least $K=\sigma/p$ disjoint pairs of indices $\{r,s\}$ such that $i \not \in \{r,s\}$ and $v_i = v_r + v_s$. It suffices to show this for any of the copies of $v_i$ in $G$, and so we may assume $1 \le i \le t$.

\noi Assume first that $G$ has no zero columns. In this case we claim that $v_i$ participates in exactly $\sum_{(j,k) \in N_i} w_j w_k$ dependencies of length three. Indeed, each pair $(j,k) \in N_i$ contributes $w_j w_k$ dependencies, obtained by taking $v_i$ together with any copy of $v_j$ and any copy of $v_k$. On the other hand, every dependency is of this form. Hence, by assumption, $\sum_{(j,k) \in N_i} w_j w_k \ge \sigma$.

\noi Next, we note that any pair $(j,k)$ in $N_i$ contributes $w_j$ disjoint pairs of indices $\{r,s\}$ such that $i \not \in \{r,s\}$ and $v_r + v_s = v_i$, obtained by making $v_r$ go over all the copies of $v_j$ in $G$ and matching each $v_r$ with a distinct copy of $v_k$. Here we use the fact that $w_j \le w_k$. Moreover, these collections of indices are disjoint for different choices of $(j,k) \in N_i$. Altogether this gives
\[
\sum_{(j,k) \in N_i} w_j \ge \frac1p \cdot \sum_{(j,k) \in N_i} w_j w_k \ge \frac{\sigma}{p}
\]
such pairs, proving the lemma in this case. For the first inequality, recall that all $w_k$ are bounded from above by $p$.

\noi If $G$ has zero columns, let $1 \leq z \leq t$ be the index with $v_z =0$. Compared to the previous case, we have $\(w_i-1\) \cdot w_z$ additional dependencies of length $3$ for $v_i$, obtained by choosing any of the extra copies of $v_i$ together with $v_i$ itself and with any copy of $v_z$. So, in this case the total number of dependencies is $\(w_i-1\) \cdot w_z + \sum_{(j,k) \in N_i} w_j w_k$, and this, by assumption, is at least~$\sigma$.

\noi On the other hand, we get $\min\{w_i-1, w_z\}$ additional disjoint pairs of indices $\{r,s\}$ such that $i \not \in \{r,s\}$ and $v_r + v_s = v_i$, by matching as many distinct copies of $v_i$ as possible (not counting $v_i$ itself) with distinct copies of $v_z$. Altogether, we get
\[
\min\{w_i-1, w_z\} + \sum_{(j,k) \in N_i} w_j ~~\ge~~ \frac1p \cdot \(\(w_i-1\) \cdot w_z +  \sum_{(j,k) \in N_i} w_j w_k\) ~~\ge~~ \frac{\sigma}{p}
\]
such pairs, proving the lemma in this case as well.

\eprf

\noi The claim of Theorem~\ref{thm:bounded_repetitions} now follows directly by substituting $K = \lceil\sigma/p\rceil$ in Corollary~\ref{cor:dimV}.

\noi We proceed with the proof of Theorem~\ref{thm:regLTC}, using the notation of Lemma~\ref{lem:bounded_repetitions}. We first note that since $C$ is a regular code, each column of $G$ has the same multiplicity $p$, implying $t = n/p$. In particular, the dimension of $C$ is at most $n/p$. Hence we may and will assume $\sigma > 4p^2$, since otherwise we are done.

\noi Next, consider the coset leader graph $\T = \{0,1\}^n / V^{\perp}$, where $V$ is the row space of $G$. By Lemmas~\ref{lem:main}~and~\ref{lem:bounded_repetitions}, the radius of $\T$ is at most $\frac{n}{\sigma/p + 1} < \frac{np}{\sigma}$. The key point to observe is that while $\T$ is an $n$-regular multigraph, the edges of $\T$ corresponding to identical columns of $G$ are parallel to each other, and hence each vertex of $\T$ has precisely $t$ distinct neighbors. Proceeding as in the proof of Corollary~\ref{cor:dimV}, we have
\[
|C| ~=~ |\T| ~\leq~ \sum_{i=0}^{\lfloor\frac{np}{\sigma}\rfloor} \binom {t}{i} ~\le~2^{tH\(\frac{np}{\sigma t}\)}.
\]
Substituting $t = n/p$, and setting $\alpha = \frac{\sigma}{p^2}$, we get
\[
\frac 1n \cdot \log_2 |C| \le \frac 1p H\(\frac {p^2}{\sigma}\) = \frac{1}{\sqrt{\sigma}} \cdot \sqrt \alpha H\(\frac {1} {\alpha}\)
\]

\noi To complete the proof, we will show that $\sqrt \alpha \cdot H(\frac{1}{\alpha}) < 2$, for all $\alpha \ge 1$. In fact,
\[
\alpha \cdot H\(\frac{1}{\alpha}\) = \log_2(\alpha) + (\alpha - 1) \log_2\(1 + \frac{1}{\alpha - 1}\) \le \frac{1}{\ln 2} \cdot \Big(\ln(\alpha) + 1\Big).
\]
Hence $\sqrt \alpha \cdot H(\frac{1}{\alpha}) \le \frac{1}{\ln 2} \cdot \frac{\ln \alpha +1 }{\sqrt \alpha}$.
It remains to observe that the function $\frac{\ln \alpha +1 }{\sqrt \alpha}$ attains its maximum of $\frac{2}{\sqrt e} < 2 \ln 2$ at $\alpha = e$.
\eprf

\noi The next example shows that Theorem~\ref{thm:bounded_repetitions} is tight, up to the $\log\(\lceil\sigma/p\rceil\)$-factor.

\xam \label{ex:counter}
\noi Let $m$ be a power of $2$, and let $k \ge \log_2 m$ be integer. Let $U$ be a linear subspace of $\{0,1\}^k$ of dimension $\log_2 m$ with minimal distance at least $3$. Let $u_1,...,u_m$ be the vectors of $U$. Let ${\bf 1}$ be the all-$1$ vector of length $k$, and let $B_i$ be the $k \times k$ matrix given by the outer product $u_i \otimes {\bf 1}$. Finally, let $I$ be the $k \times k$ identity matrix.

\noi

\noi Let $G$ be the following $k \times n$ matrix with $n = 2km$. The first $km$ columns of $G$ are formed by $m$ square blocks $I + B_i$, for $i = 1,...,m$. The remaining $km$ columns are formed by the blocks $B_1,...,B_m$.

\noi Clearly the rows of $G$ are linearly independent, and therefore the dimension of the code $C$ it generates is $k$. By construction, for $G$ holds $p = k$ and $\sigma = km$ (since $U$ is a subspace). Hence we have
\[
dim(C) = k = \frac{n}{2m} = \frac{n}{2 \sigma/p}.
\]
\exam

\bibliographystyle{alpha}
\bibliography{../eran}

\newcommand{\etalchar}[1]{$^{#1}$}
\begin{thebibliography}{KKRT15}

\bibitem[BHL{\etalchar{+}}15]{bauer2015li}
Frank Bauer, Paul Horn, Yong Lin, Gabor Lippner, Dan Mangoubi, and Shing-Tung
  Yau.
\newblock Li-{Y}au inequality on graphs.
\newblock {\em J. Differential Geom.}, 99(3):359--405, 2015.

\bibitem[BJL12]{bauer2012ollivier}
Frank Bauer, J{\"u}rgen Jost, and Shiping Liu.
\newblock Ollivier-{R}icci curvature and the spectrum of the normalized graph
  {L}aplace operator.
\newblock {\em Math. Res. Lett.}, 19(6):1185--1205, 2012.

\bibitem[BSV12]{ben2012towards}
Eli Ben-Sasson and Michael Viderman.
\newblock Towards lower bounds on locally testable codes via density arguments.
\newblock {\em Comput. Complexity}, 21(2):267--309, 2012.

\bibitem[Cha96]{charney1994metric}
Ruth Charney.
\newblock Metric geometry: connections with combinatorics.
\newblock In {\em Formal power series and algebraic combinatorics ({N}ew
  {B}runswick, {NJ}, 1994)}, volume~24 of {\em DIMACS Ser. Discrete Math.
  Theoret. Comput. Sci.}, pages 55--69. Amer. Math. Soc., Providence, RI, 1996.

\bibitem[CHLL97]{cohen1997covering}
G{\'e}rard Cohen, Iiro Honkala, Simon Litsyn, and Antoine Lobstein.
\newblock {\em Covering codes}, volume~54 of {\em North-Holland Mathematical
  Library}.
\newblock North-Holland Publishing Co., Amsterdam, 1997.

\bibitem[DK11]{dinur2011dense}
Irit Dinur and Tali Kaufman.
\newblock Dense locally testable codes cannot have constant rate and distance.
\newblock In {\em Approximation, randomization, and combinatorial
  optimization}, volume 6845 of {\em Lecture Notes in Comput. Sci.}, pages
  507--518. Springer, Heidelberg, 2011.

\bibitem[DSW14]{dvir2014breaking}
Zeev Dvir, Shubhangi Saraf, and Avi Wigderson.
\newblock Breaking the quadratic barrier for 3-lcc's over the reals.
\newblock In {\em Proceedings of the forty-sixth annual ACM symposium on Theory
  of computing}, pages 784--793. ACM, 2014.

\bibitem[FT05]{friedman2006generalized}
Joel Friedman and Jean-Pierre Tillich.
\newblock Generalized {A}lon-{B}oppana theorems and error-correcting codes.
\newblock {\em SIAM J. Discrete Math.}, 19(3):700--718 (electronic), 2005.

\bibitem[GKST06]{goldreich2006lower}
Oded Goldreich, Howard Karloff, Leonard~J. Schulman, and Luca Trevisan.
\newblock Lower bounds for linear locally decodable codes and private
  information retrieval.
\newblock {\em Comput. Complexity}, 15(3):263--296, 2006.

\bibitem[KdW04]{kerenidis2004exponential}
Iordanis Kerenidis and Ronald de~Wolf.
\newblock Exponential lower bound for 2-query locally decodable codes via a
  quantum argument.
\newblock {\em J. Comput. System Sci.}, 69(3):395--420, 2004.

\bibitem[KKRT15]{klartag2015discrete}
Bo'az Klartag, Gady Kozma, Peter Ralli, and Prasad Tetali.
\newblock Discrete curvature and abelian groups.
\newblock {\em arXiv preprint arXiv:1501.00516}, 2015.

\bibitem[KT00]{katz2000efficiency}
Jonathan Katz and Luca Trevisan.
\newblock On the efficiency of local decoding procedures for error-correcting
  codes.
\newblock In {\em Proceedings of the {T}hirty-{S}econd {A}nnual {ACM}
  {S}ymposium on {T}heory of {C}omputing}, pages 80--86 (electronic), New York,
  2000. ACM.

\bibitem[KV10]{kaufman2010locally}
Tali Kaufman and Michael Viderman.
\newblock Locally testable vs. locally decodable codes.
\newblock In {\em Approximation, randomization, and combinatorial
  optimization}, volume 6302 of {\em Lecture Notes in Comput. Sci.}, pages
  670--682. Springer, Berlin, 2010.

\bibitem[LPW09]{levin2009markov}
David~A. Levin, Yuval Peres, and Elizabeth~L. Wilmer.
\newblock {\em Markov chains and mixing times}.
\newblock American Mathematical Society, Providence, RI, 2009.
\newblock With a chapter by James G. Propp and David B. Wilson.

\bibitem[LY10]{lin2010ricci}
Yong Lin and Shing-Tung Yau.
\newblock Ricci curvature and eigenvalue estimate on locally finite graphs.
\newblock {\em Math. Res. Lett.}, 17(2):343--356, 2010.

\bibitem[MS77]{macwilliams1977theory}
F.~J. MacWilliams and N.~J.~A. Sloane.
\newblock {\em The theory of error-correcting codes}.
\newblock North-Holland Publishing Co., Amsterdam-New York-Oxford, 1977.
\newblock North-Holland Mathematical Library, Vol. 16.

\bibitem[Oll09]{ollivier2009ricci}
Yann Ollivier.
\newblock Ricci curvature of {M}arkov chains on metric spaces.
\newblock {\em J. Funct. Anal.}, 256(3):810--864, 2009.

\bibitem[OV12]{ollivier2012curved}
Y.~Ollivier and C.~Villani.
\newblock A curved {B}runn-{M}inkowski inequality on the discrete hypercube,
  or: what is the {R}icci curvature of the discrete hypercube?
\newblock {\em SIAM J. Discrete Math.}, 26(3):983--996, 2012.

\bibitem[Pet11]{petrunin2011alexandrov}
Anton Petrunin.
\newblock Alexandrov meets {L}ott-{V}illani-{S}turm.
\newblock {\em M\"unster J. Math.}, 4:53--64, 2011.

\bibitem[vL99]{van1999introduction}
J.~H. van Lint.
\newblock {\em Introduction to coding theory}, volume~86 of {\em Graduate Texts
  in Mathematics}.
\newblock Springer-Verlag, Berlin, third edition, 1999.

\bibitem[Woo07]{woodruff2007new}
D.~Woodruff.
\newblock New lower bounds for general locally decodable codes.
\newblock In {\em Electronic Colloquium on Computational Complexity (ECCC)},
  volume~14, 2007.

\bibitem[Woo12]{woodruff2012quadratic}
David~P. Woodruff.
\newblock A quadratic lower bound for three-query linear locally decodable
  codes over any field.
\newblock {\em J. Comput. Sci. Tech.}, 27(4):678--686, 2012.

\end{thebibliography}

\end{document}